\documentclass{article}
\usepackage{graphicx}
\usepackage{amsmath}
\usepackage{amsfonts}
\allowdisplaybreaks[4]
\newtheorem{theorem}{Theorem}

\title{A closed solution to a special polynomial trinomial equation and semi-analytical roots for a general algebraic equation}
\author{Rong Qiang Wei\thanks{
                  College of Earth and Planet Sciences, University of Chinese Academy of Sciences, Beijing, PRC, 100049.
                  e\_mail: wrq1973@ucas.edu.cn}}
\date{}
\begin{document}
\maketitle

\begin{abstract}

   We suggest a closed solution for the roots of polynomial trinomial algebraic equation $$z^n+xz^{n-1}-1=0$$ with an appropriate $x$. This solution is a minor modification to the work of Mikhalkin (Mikhalkin E N, 2006. On solving general algebraic equations by integrals of elementary functions, Siberian Mathematical Jounral, 47(2), 301-306). This modification, together with Mikhalkin's integral formula, provides a relatively simple analytical expression for the solution to a general algebraic equation when the polynomial coefficients are over the corresponding convergent domain. Numerical examples show that this expression can be another alternative to finding numerically the roots of a general polynomial algebraic equation when the integral involved exists and is calculated correctly.     
      
\end{abstract}

{\hspace{2.2em}\small Keywords:}

{\hspace{2.2em}\tiny polynomial algebraic equation, closed solution, integral expression, analytic continuation}

\section{Introduction}

The identification of roots for the polynomial algebraic equation is an important problem in many scientific and technical applications. The related studies have a long history of thousand years, which can be traced from ancient Babylonia. For the quadratic equations, 

\begin{equation}\label{eq1}
z^2+x_1z+x_2=0
\end{equation}
where $z, x_i\in \mathbb{C}$, there exists a beautiful root formula (Eq. (\ref{eq2})), 

\begin{equation}\label{eq2}
z_{_{1,2}}=\frac{-x_1\pm\sqrt{x_1^2-4x_2}}{2}
\end{equation}

\begin{figure}
	\begin{center}
		\includegraphics[scale=0.6]{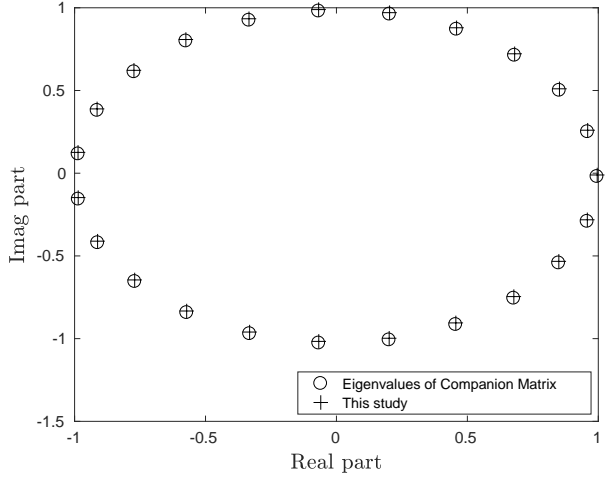}
	\end{center}
	\caption{Roots of $z^{23}+0.5iz^{22}-1=0$ from Eq. (\ref{wrq_intg2}) and those from the eigenvalues of corresponding companion matrix.}
	\label{fig_z23}
\end{figure}

For the cubic equations, 

\begin{equation}\label{eq3}
z^3+x_1z^2+x_2z+x_3=0
\end{equation}
there also exists a root formula (Zeidler et al., 2004), but it is complicated. Firstly, let $y=z+\frac{x_1}{3}$, and Eq. (\ref{eq3}) can be transformed into,

\begin{equation}\label{eq4}
y^3+3py+2q=0
\end{equation}
where $3p=x_2-\frac{x_1^2}{3}$,$2q=\frac{2x_1^3}{27}-\frac{x_1x_2}{3}+x_3$.
 
 Secondly Eq. (\ref{eq4}) has three roots ({\it Cardano's formula}),
 
\begin{equation}\label{eq5}
\left\{ \begin{array}{l}
	y_1=u_++u_-\\
	y_2=\rho _+u_++\rho _-u_-\\
	y_3=\rho _-u_++\rho _+u_-\\
\end{array} \right. 
\end{equation} 
where 
$$
u_{\pm}:= \sqrt[3]{-q\pm \sqrt{p^3+q^2}}
$$
and
$$
\rho _{\pm}:= \frac{1}{2}\left( -1\pm i\sqrt{3} \right) 
$$

Finally, the roots of Eq. (\ref{eq3}) are,

\begin{equation}\label{eq6}
z_{i}=y_i-\frac{x_1}{3}
\end{equation} 
where $i=1,2,3$.

For the biquadratic equations,

\begin{equation}\label{eq7}
z^4+x_1z^3+x_2z^2+x_3z+x_4=0
\end{equation}
there still exists a root formula (Zeidler et al., 2004), but it is too complicated to be shown explicitly. Again firstly Eq. (\ref{eq7}) can be brought into a normal form when we let $y=z+\frac{x_1}{4}$,

\begin{equation}\label{eq8}
y^4+py^2+qy+r=0
\end{equation}

\begin{figure}
	\begin{center}
		\includegraphics[scale=0.75]{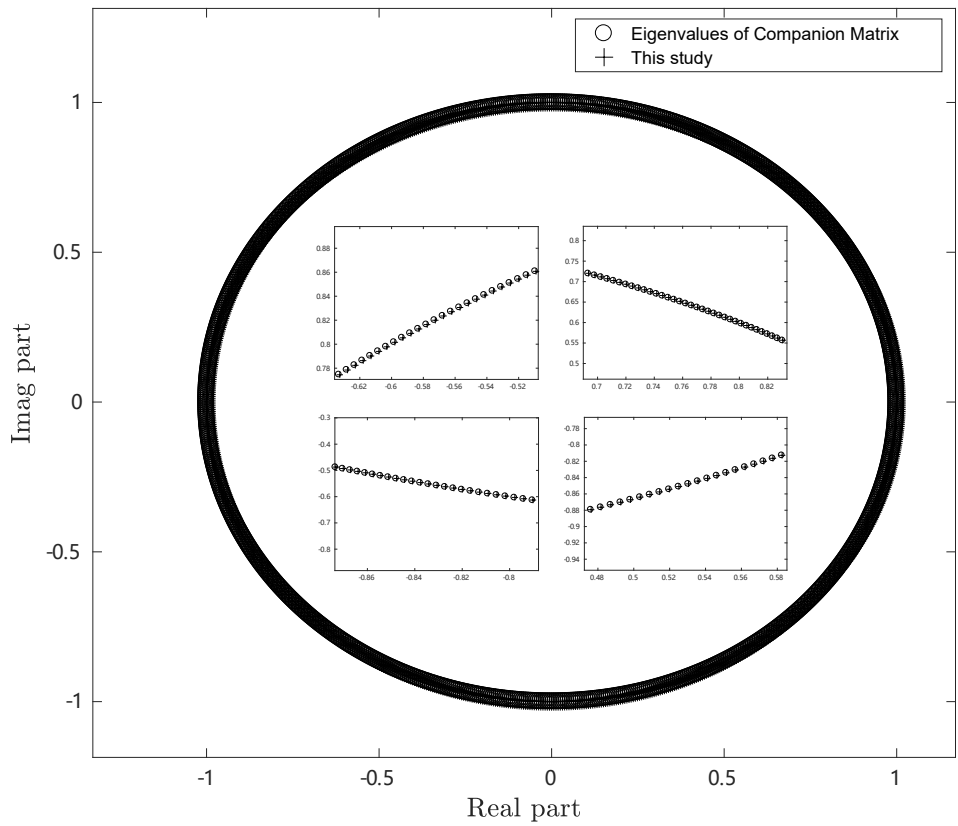}
	\end{center}
	\caption{Roots of $z^{1000}+(0.5-0.37i)z^{999}-1=0$ from Eq. (\ref{wrq_intg2}) and those from the eigenvalues of corresponding companion matrix. Four subfigures show the detailed comparison for partial roots of these two types.}
	\label{fig_z1000}
\end{figure}

\begin{figure}\label{fig_z3}
	\begin{center}
		\includegraphics[scale=0.6]{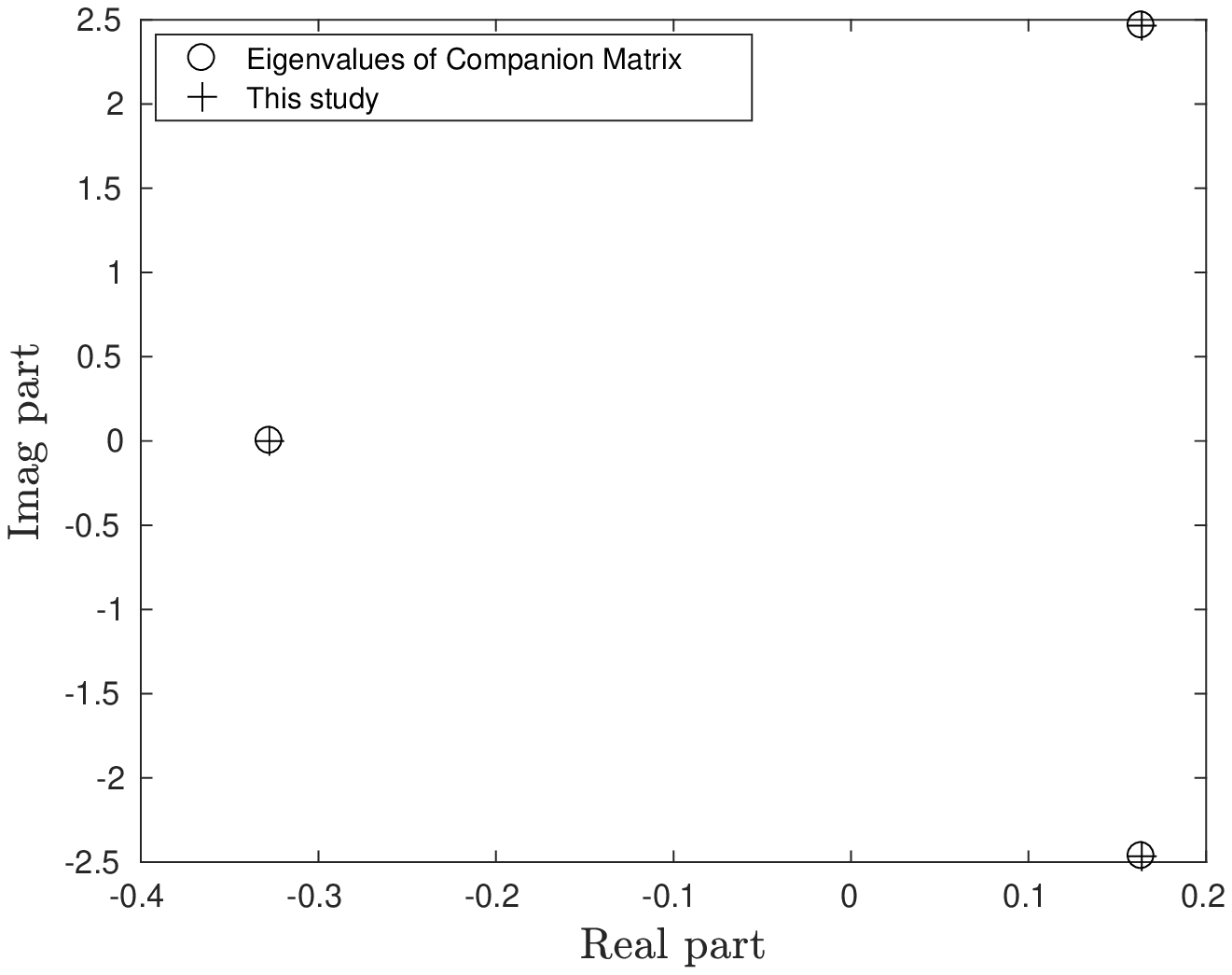}
	\end{center}
	\caption{Roots of $z^3+6z+2=0$ from Eq. (\ref{Mikhalkin_intg1}) and those from the eigenvalues of corresponding companion matrix.}
\end{figure}

Secondly, the roots of Eq. (\ref{eq8}) are,

\begin{equation}\label{eq9}
\left\{ \begin{array}{l}
	2y_1=u+v+w\\
	2y_2=u-v+w\\
	2y_3=-u+v+w\\
	2y_4=-u-v-w\\
\end{array} \right. 
\end{equation}
where $u$,$v$,$w$ are solutions of the equations $u^2=\alpha$, $v^2=\beta$ and $w^2=\gamma$ (where $uvw=q$ is required), respectively. And $\alpha,\beta,\gamma$ are the roots of the cubic equation ({\it cubic resolvent}) in the following,

$$
z^3+2pz^2+\left( p^2-4r \right) z-q^2=0
$$
 
Finally, the roots of Eq. (\ref{eq7}) are,

\begin{equation}\label{eq10}
z_{i}=y_i-\frac{x_1}{4}
\end{equation} 
where $i=1,2,3,4$.

\begin{figure}
	\begin{center}
		\includegraphics[scale=0.6]{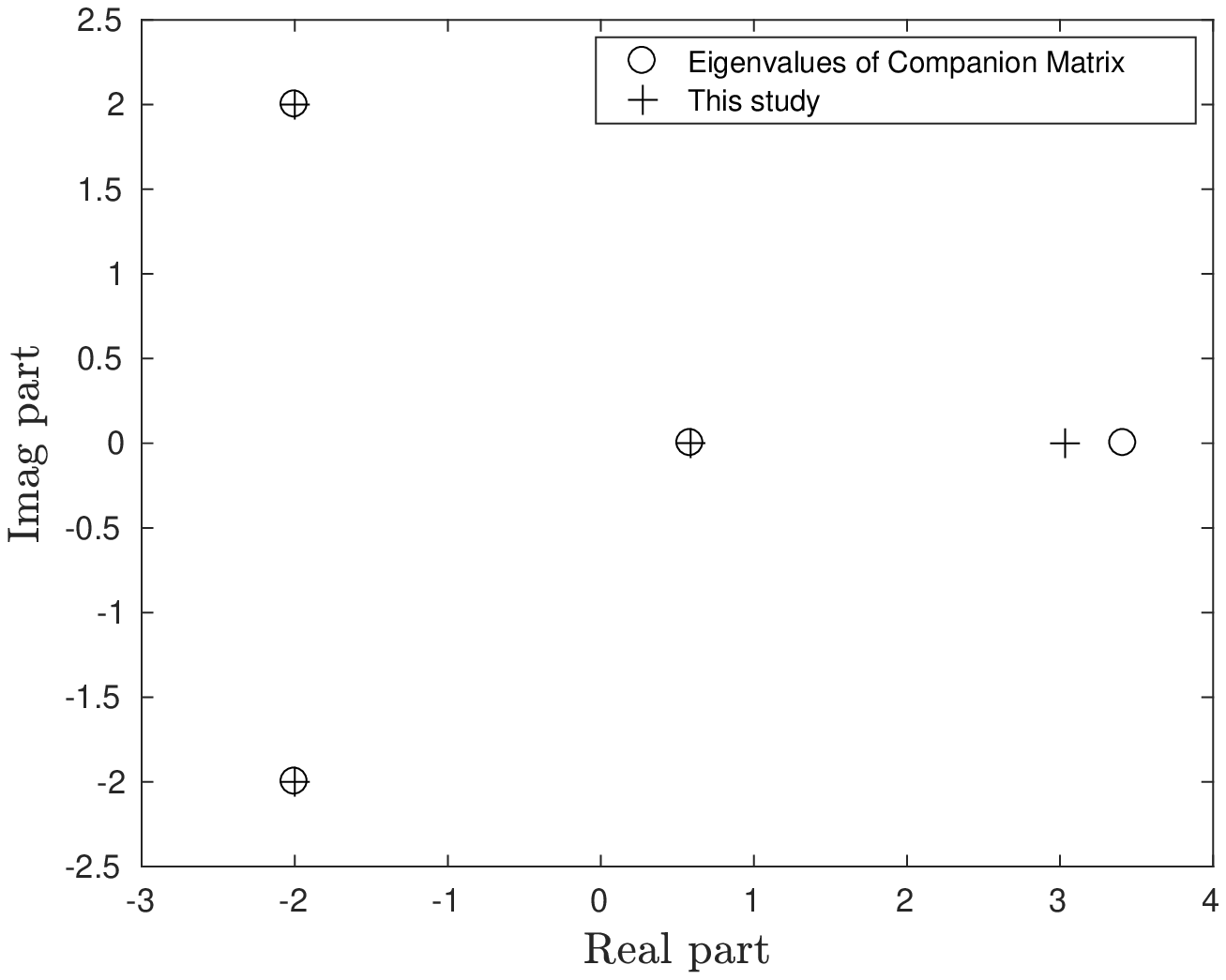}
	\end{center}
	\caption{Roots of $z^4-6z^2-24z+16=0$ from Eq. (\ref{Mikhalkin_intg1}) and those from the eigenvalues of corresponding companion matrix.}
	\label{fig_z4}
\end{figure}

It can be seen from the above that the roots of polynomial algebraic equations of $n^{th}$ degree with $n=2,3,4$ can be expressed explicitly or implicitly only by the radicals. However, this does not hold for the general algebraic equations with $n\geq 5$, which has been proven by the great Abel and Galois. However, it is shown from studies later that there exist other ways to obtain the roots of  these equations. The most common way is to solve them numerically. There perhaps have been hundreds if not thousands of different methods devised for numerical solution to the roots of polynomials (McNamee and Pan, 2012), and many are described in McNamee (2007) and McNamee and Pan (2013). Recently, new studies in this area are still emerging (eg., Raffalli, 2022). Among them, one of the most notable method is computing the eigenvalues of the corresponding companion matrix of polynomials with QR algorithm, by which it is possible to obtain all the roots both quickly and accurately.

On the other hand, people have never stopped trying to solve the algebraic equations of $n^{th}$ degree with $n\geq 5$ using analytical methods, which is meaningful not only in theory but also in numerical calculations. Among them, trinomial equations ($z^n+xz^m-1=0, 0<m<n$) are undoubtedly special, but are no less general. For the trinomial equation of the form $z^n+z=x$, Glasser (2000) and Ritelli and Spaletta (2021) solved it in hypergeometric way with or without the Lagrange's implicit functions theorem, respectively. In case of general $n^{th}$ degree equation, Umemura (1984) proved a formula for it using the theory of special function, which can be stated as the follows,

\begin{theorem}\label{Theo1}
    Let $$
     P\left( x \right) =a_0x^n+a_1x^{n-1}+\cdots +a_n=0, a_0\neq 0, a_i\in\mathbb{C} (0\leq i\leq n), 
    $$
be an algebraic equation irreducible over a certain subfield of $\mathbb{C}$, then a root of it is given by 
\begin{equation}\label{eqTheo1}
\begin{array}{c}
	\left[ \theta \left( \begin{matrix}
	\frac{1}{2}&		0&		\cdots&		0\\
	0&		0&		\cdots&		0\\
\end{matrix} \right) \left( \varOmega \right) ^4\theta \left( \begin{matrix}
	\frac{1}{2}&		\frac{1}{2}&		\cdots&		0\\
	0&		0&		\cdots&		0\\
\end{matrix} \right) \left( \varOmega \right) ^4 \right.\\
	+\theta \left( \begin{matrix}
	0&		0&		\cdots&		0\\
	0&		0&		\cdots&		0\\
\end{matrix} \right) \left( \varOmega \right) ^4\theta \left( \begin{matrix}
	0&		\frac{1}{2}&		\cdots&		0\\
	0&		0&		\cdots&		0\\
\end{matrix} \right) \left( \varOmega \right) ^4\\
	\left. -\theta \left( \begin{matrix}
	0&		0&		\cdots&		0\\
	\frac{1}{2}&		0&		\cdots&		0\\
\end{matrix} \right) \left( \varOmega \right) ^4\theta \left( \begin{matrix}
	0&		\frac{1}{2}&		\cdots&		0\\
	\frac{1}{2}&		0&		\cdots&		0\\
\end{matrix} \right) \left( \varOmega \right) ^4 \right]\\
	/\left[ 2\theta \left( \begin{matrix}
	\frac{1}{2}&		0&		\cdots&		0\\
	0&		0&		\cdots&		0\\
\end{matrix} \right) \left( \varOmega \right) ^4\theta \left( \begin{matrix}
	\frac{1}{2}&		\frac{1}{2}&		\cdots&		0\\
	0&		0&		\cdots&		0\\
\end{matrix} \right) \left( \varOmega \right) ^4 \right]\\
\end{array}
 \end{equation}
  where $$
\theta \left( \begin{array}{c}
	m_1\\
	m_2\\
\end{array} \right) \left( z,\tau \right) =\sum_{\xi \in \mathbb{Z}^g}{e^{\left( \frac{1}{2}\left( \xi +m_1 \right) \tau ^T\left( \xi +m_1 \right) +\left( \xi +m_1 \right) ^T\left( z+m_2 \right) \right)}}
$$
is the theta function defined for row vector $m_1,m_2\in\mathbb{R}^g$, $z\in\mathbb{C}^g$ and a symmetric $g\times g$ matrix $\tau$ with positive definite imaginary part, and here $e(x)=\exp(2\pi i x)$; $\Omega$ is the period matrix of the hyperelliptic curve $y^2=F(x)$ ($F(x)$ has $2g+1$ simple roots), and $F(x)=x(x-1)P(x)$ if the degree of $P(x)$ is odd, and $F(x)=x(x-1)(x-2)P(x)$, otherwise.
\end{theorem}

\begin{figure}
	\begin{center}
		\includegraphics[scale=0.6]{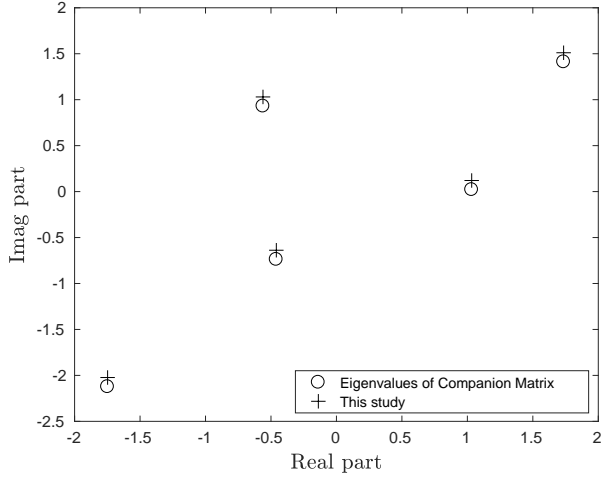}
	\end{center}
	\caption{Roots of $z^5+0.5iz^4-6iz^3-2.4z^2+z+6i=0$ from Eq. (\ref{Mikhalkin_intg1}) and those from the eigenvalues of corresponding companion matrix.}
	\label{fig_z5}
\end{figure}

  However, it can be seen that  {\bf Theorem \ref{Theo1}}, ie., Eq. (\ref{eqTheo1}) is too complicated and technical to be practically useful, especially it may be not easy to obtain the period matrix $\varOmega$. 
  
  In what follows we will introduce another analytical solution by Mellin (1921) to a general algebraic equation, which is simpler than Eq. (\ref{eqTheo1}) although it has a limited domain of convergence. Further, we will introduce a simplified version of this solution using the integrals of elementary functions by Mikhalkin (2006), and a minor modification of ours to Mikhalkin (2006) for finding the roots of the trinomial algebraic equation $z^n+xz^{n-1}-1=0$. Finally, we will show some numerical examples based on Mikhalkin's integral formula and our minor modification.         

\section{Mellin's solution in integral and series}\label{Mellin1921}

Mellin (1921) obtained an integral formula to compute the principal root (branch)
of the polynomial equation $z^n + x_1z^{n_1} + \cdots + x_pz^{n_p}-1=0$ using hypergeometric functions of its coefficients $x_1, ..., x_p$ with the condition $z(0,0,...0)=1$. This formula is as the follows when $\vert {\rm{arg}}(x_m) \vert < n_p\pi/{2n}, m=1,2,...p$,

 \begin{align}\label{Mellin1921_intg1}
  \begin{aligned}
z\left( x_1,x_2,\cdots ,x_p \right) =&		\frac{1}{\left( 2\pi i \right) ^p}\int_{\gamma +i\mathbb{R} ^p}{\frac{\frac{1}{n}\mathbf{\Gamma }\left( \frac{1}{n}-\frac{n_1}{n}z_1-\frac{n_2}{n}z_2-\cdots -\frac{n_p}{n}z_p \right)}{\mathbf{\Gamma }\left( \frac{1}{n}+\frac{n_{1}^{'}}{n}z_1+\frac{n_{2}^{'}}{n}z_2+\cdots +\frac{n_{p}^{'}}{n}z_p+1 \right)}}\\
	&		\mathbf{\Gamma }\left( z_1 \right) \mathbf{\Gamma }\left( z_2 \right) \cdots \mathbf{\Gamma }\left( z_p \right) x_{1}^{-z_1}x_{2}^{-z_2}\cdots x_{p}^{-z_p}\mathrm{d}z_1\mathrm{d}z_2\cdots \mathrm{d}z_p\\
	&		\\
\end{aligned}
 \end{align}
where $\gamma$ is a point in the polyhedron, $$
\left\{ u\in \mathbb{R} ^p:u_1>0,\cdots u_p>0,n_1u_1+n_2u_2+\cdots n_pu_p<1 \right\}$$
and $n_{m}^{'}=n-n_m,m=1,2,\cdots,p$.

The other $n-1$ roots (branches) are obtained from $z\left( x_1,x_2,\cdots ,x_p \right)$ by the formula,

\begin{equation}\label{other_roots}
z_j\left( x_1,x_2,\cdots ,x_p \right) =\varepsilon ^jz\left( x_1\varepsilon ^{jn_1},x_2\varepsilon ^{jn_2},\cdots ,x_p\varepsilon ^{jn_p} \right) ,     j=1,2,\cdots ,n-1.
\end{equation}
where $\varepsilon=\exp (\frac{2\pi i}{n})$ is a primitive $n$th root of unity.

This paper is written in French in a terse style, and an relatively elementary explanation or an expository article can be Lawton (2021). 

Mellin (1921) also presented the formula for any positive powers of the principal root, ie., for $\mu >0$,

\begin{align}\label{Mellin1921_intg2}
  \begin{aligned}
z^{\mu}\left( x_1,x_2,\cdots ,x_p \right) =&		\frac{1}{\left( 2\pi i \right) ^p}\int_{\gamma +i\mathbb{R} ^p}{\frac{\frac{\mu}{n}\mathbf{\Gamma }\left( \frac{\mu}{n}-\frac{n_1}{n}z_1-\frac{n_2}{n}z_2-\cdots -\frac{n_p}{n}z_p \right)}{\mathbf{\Gamma }\left( \frac{\mu}{n}+\frac{n_{1}^{'}}{n}z_1+\frac{n_{2}^{'}}{n}z_2+\cdots +\frac{n_{p}^{'}}{n}z_p+1 \right)}}\\
	&		\mathbf{\Gamma }\left( z_1 \right) \mathbf{\Gamma }\left( z_2 \right) \cdots \mathbf{\Gamma }\left( z_p \right) x_{1}^{-z_1}x_{2}^{-z_2}\cdots x_{p}^{-z_p}\mathrm{d}z_1\mathrm{d}z_2\cdots \mathrm{d}z_p\\
\end{aligned}
 \end{align}  
where $\gamma$ is a point in the polyhedron, $$
\left\{ u\in \mathbb{R} ^p:u_1>0,\cdots u_p>0,n_1u_1+n_2u_2+\cdots n_pu_p<\mu \right\}$$

Calculating the integral (\ref{Mellin1921_intg1}) as the sum of residues at the poles of $z=-k, k=0,1,2,....$ of Euler gamma function $\Gamma(z_p), p=1,2,...$, Mellin (1921) presented the representation of $z\left( x_1,x_2,\cdots ,x_p \right)$ in the form of a hypergeometric series,

 \begin{align}\label{Mellin1921_series1}
  \begin{aligned}
z\left( x_1,x_2,\cdots ,x_p \right) =&
\frac{1}{n}\sum_{\left| k \right|\geq 0}{\frac{\left( -1 \right) ^{\left| k \right|}\mathbf{\Gamma }\left( \frac{1}{n}+\frac{n_1}{n}k_1+\frac{n_2}{n}k_2+\cdots +\frac{n_p}{n}k_p \right) x_{1}^{k_1}x_{2}^{k_2}\cdots x_{p}^{k_p}}{k_1!k_2!\cdots k_p!\mathbf{\Gamma }\left( \frac{1}{n}-\frac{n_{1}^{'}}{n}k_1-\frac{n_{2}^{'}}{n}k_2-\cdots -\frac{n_{p}^{'}}{n}k_p+1 \right)}}\\
\end{aligned}
 \end{align}  
where $\vert k\vert =k_1+k_2+\cdots+k_p; k_1, k_2, \cdots k_p\geq 0$.

Correspondingly,

\begin{align}\label{Mellin1921_series2}
  \begin{aligned}
z^\mu\left( x_1,x_2,\cdots ,x_p \right) =&
\frac{\mu}{n}\sum_{\left| k \right|\geq 0}{\frac{\left( -1 \right) ^{\left| k \right|}\mathbf{\Gamma }\left( \frac{\mu}{n}+\frac{n_1}{n}k_1+\frac{n_2}{n}k_2+\cdots +\frac{n_p}{n}k_p \right) x_{1}^{k_1}x_{2}^{k_2}\cdots x_{p}^{k_p}}{k_1!k_2!\cdots k_p!\mathbf{\Gamma }\left( \frac{\mu}{n}-\frac{n_{1}^{'}}{n}k_1-\frac{n_{2}^{'}}{n}k_2-\cdots -\frac{n_{p}^{'}}{n}k_p+1 \right)}}\\
\end{aligned}
 \end{align} 

From the view point of theory, we can solve analytically the general algebraic equations of $n^{th}$ degree with $n\geq 5$ through Mellin's solution in integral and/or series when $x_i,i=1, 2,\cdots p$ are over the convergent domain of these integral and/or series. 

\begin{figure}
	\begin{center}
		\includegraphics[scale=0.6]{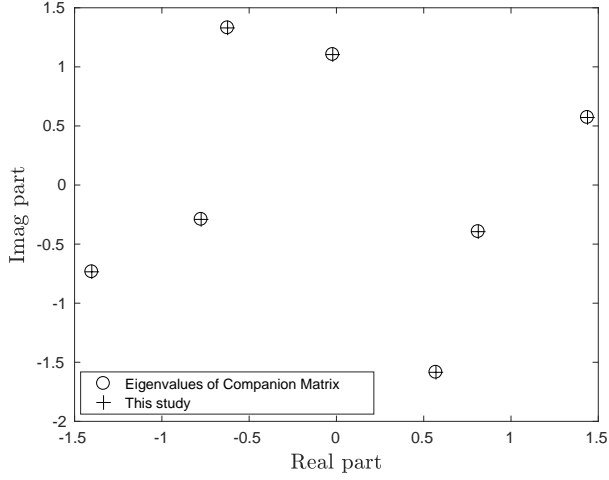}
	\end{center}
	\caption{Roots of $z^7+0.3z^5+0.5iz^4-6iz^3-2.4z^2+z+5=0$ from Eq. (\ref{Mikhalkin_intg1}) and those from the eigenvalues of corresponding companion matrix.}
	\label{fig_z7}
\end{figure}

\section{Mikhalkin's integral formula}

However, it is clear that the integrand in the integral (\ref{Mellin1921_intg1}) is a transcendental function and the integration domain is unbounded, while the solution (\ref{Mellin1921_series1}) is a multiple series, both of which are difficult to be calculated in a normal way. Mikhalkin (2006) proposed an elegant and simple integral formula based on series (\ref{Mellin1921_series1}), in which only the 1D integration of elementary functions over $[0,1]$ is involved. This formula, which has a broader convergent domain than that of the integral (\ref{Mellin1921_intg1}), reads,

\begin{align}\label{Mikhalkin_intg1}
  \begin{aligned}
z\left( x_1,x_2,\cdots ,x_p \right) =&1+\frac{1}{\left( 2\pi i \right) n}\int_0^1{t^{\frac{1-n}{n}}}\left( 1-t \right) ^{-\frac{1+n}{n}}\left[ e^{\frac{\pi i}{n}}\log \left( 1+\sum_{k=1}^p{y_ke^{-\frac{n_{k}^{'}}{n}\pi i}} \right) \right. \\
& -\left. e^{-\frac{\pi i}{n}}\log \left( 1+\sum_{k=1}^p{y_ke^{\frac{n_{k}^{'}}{n}\pi i}} \right) \right] \mathrm{d}t \\
\end{aligned}
 \end{align}  
where $n_{k}^{'}=n-n_k$, $y_k=x_kt^{\frac{n_k}{n}}\left( 1-t \right) ^{\frac{n_{k}^{'}}{n}}$.

Series (\ref{Mellin1921_series1}) is convergent only for the small $\vert x_i\vert, i=1,2,...p$, while Eq. (\ref{Mikhalkin_intg1}) is valid for any (complex) $x_i$, except $\varSigma _-\cup \varSigma _+ $ where $$\varSigma _-=\bigcup_{t\in \left[ 0;1 \right]}{\{ \sum_{k=1}^p{y_ke^{-\frac{n_{k}^{'}}{n}\pi i}}+1=0 \}}$$ and $$\varSigma _+=\bigcup_{t\in \left[ 0;1 \right]}{\{ \sum_{k=1}^p{y_ke^{\frac{n_{k}^{'}}{n}\pi i}}+1=0\}}$$

 Therefore, $z\left( x_1,x_2,\cdots ,x_p \right)$ which comes from the Integral (\ref{Mikhalkin_intg1}) is holomorphic and has a univalent continuation from a neighborhood of the origin to the domain 
$\mathbb{C} ^p\backslash \left( \varSigma _-\cup \varSigma _+ \right)$. 

\begin{figure}
	\begin{center}
		\includegraphics[scale=0.6]{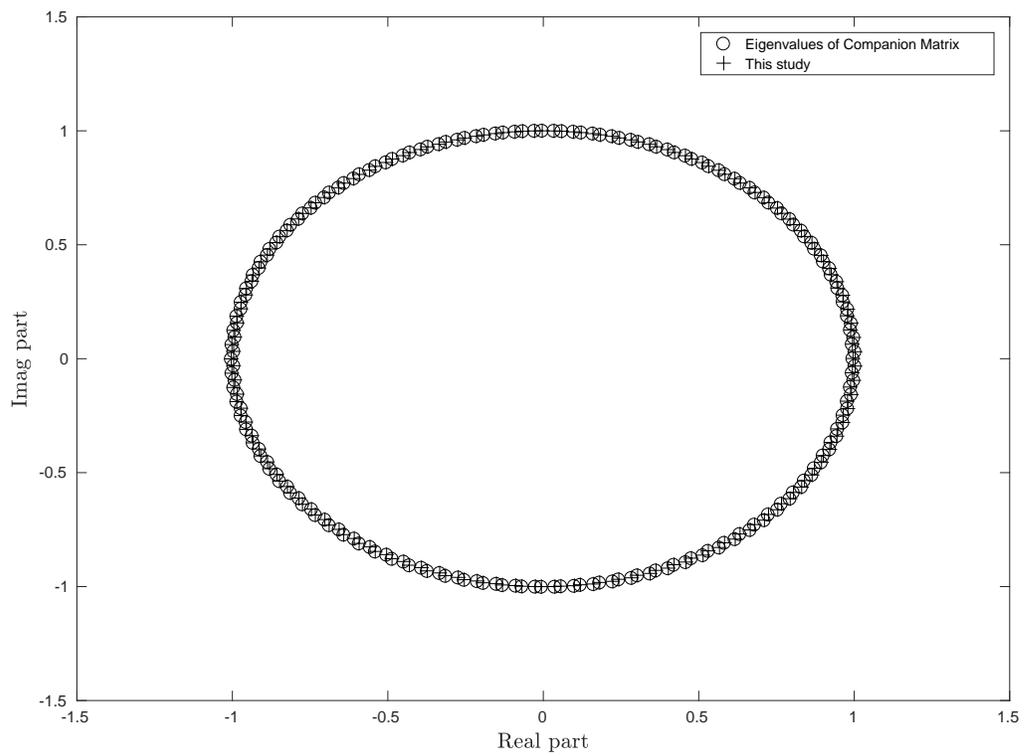}
	\end{center}
	\caption{Roots of $z^{200}+0.5z^{199}+0.8z^{99}-1=0$ from Eq. (\ref{Mikhalkin_intg1}) and those from the eigenvalues of corresponding companion matrix.}
	\label{fig_z200}
\end{figure}  

\section{A minor modification to Mikhalkin's integral formula}

Wrong roots will be obtained when the formula (\ref{Mikhalkin_intg1}) is applied to the trinomial algebraic equation $z^n+xz^{n-1}-1=0$. Except that $x$ is in a non-convergent domain of this integral, the reason is from the Eq. (\ref{f_beta}) which is used in the deriving Integral formula (\ref{Mikhalkin_intg1}), 

\begin{align}\label{f_beta}
  \begin{aligned}
\frac{\Gamma \left( x \right) \Gamma \left( y \right)}{\Gamma \left( x+y \right)}=\mathrm{B}\left( x,y \right) =\int_0^1{t^{x-1}}\left( 1-t \right) ^{y-1}\mathrm{d}t
\end{aligned}
 \end{align}  

The right hand of Eq. (\ref{f_beta}) requires ${\rm Re\ } x>0, {\rm Re\ } y>0$, although there is no such a requirement when we evaluate the integral of $\int_0^1{t^{x-1}}\left( 1-t \right) ^{y-1}\mathrm{d}t$ from $\frac{\Gamma \left( x \right) \Gamma \left( y \right)}{\Gamma \left( x+y \right)}$. To meet this requirement, we re-derive formula (\ref{Mikhalkin_intg1}) from series solution (\ref{Mellin1921_series2}) other than (\ref{Mellin1921_series1}). We obtain, 

\begin{align}\label{wrq_intg1}
  \begin{aligned}
z^\mu\left( x_1,x_2,\cdots ,x_p \right) =&1+\frac{\mu}{\left( 2\pi i \right) n}\int_0^1{t^{\frac{\mu}{n}-1}}\left( 1-t \right) ^{-\frac{\mu}{n}-1}\left[ e^{\frac{\mu}{n}\pi i}\log \left( 1+\sum_{k=1}^p{y_ke^{-\frac{n_{k}^{'}}{n}\pi i}} \right) \right. \\
& -\left. e^{-\frac{\mu}{n}\pi i}\log \left( 1+\sum_{k=1}^p{y_ke^{\frac{n_{k}^{'}}{n}\pi i}} \right) \right] \mathrm{d}t \\
\end{aligned}
 \end{align}  
 
 For this case, we found that in Beta function  $y=-\frac{\mu}{n}+\frac{n_{1}^{'}}{n}k_1+\frac{n_{2}^{'}}{n}k_2+\cdots +\frac{n_{p}^{'}}{n}k_p$; If $\mu=1$, $y=0$ for algebraic equation $z^n+xz^{n-1}-1=0$ when $k_1=1$, which is not the requirement of ${\rm Re\ } y>0$  in Eq. (\ref{f_beta}); But if $0<\mu <1$, $y>0$. For convenience we let $\mu = \frac{1}{2}$ here. So now the principal root (branch) for this trinomial equation is,

\begin{align}\label{wrq_intg2}
  \begin{aligned}
z(x) =&\left\{1+\frac{1}{\left( 4\pi i \right) n}\int_0^1{t^{\frac{1}{2n}-1}}\left( 1-t \right) ^{-\frac{1}{2n}-1}\left[ e^{\frac{\pi i}{2n}}\log \left( 1+xt^{\frac{n-1}{n}}\left( 1-t \right) ^{\frac{1}{n}}e^{-\frac{\pi i}{n}} \right) \right.\right. \\
& -\left. \left. e^{-\frac{\pi i}{2n}}\log \left( 1+xt^{\frac{n-1}{n}}\left( 1-t \right) ^{\frac{1}{n}}e^{\frac{\pi i}{n}} \right) \right] \mathrm{d}t\right\}^2 \\
\end{aligned}
 \end{align}  
where the convergent domain of $x$ is $\mathbb{C} ^p\backslash \left( \varSigma' _-\cup \varSigma' _+ \right)$, and $$\varSigma' _-=\bigcup_{t\in \left[ 0;1 \right]}{\{ 1+xt^{\frac{n-1}{n}}\left( 1-t \right) ^{\frac{1}{n}}e^{-\frac{\pi i}{n}}=0 \}}$$ and $$\varSigma' _+=\bigcup_{t\in \left[ 0;1 \right]}{\{ 1+xt^{\frac{n-1}{n}}\left( 1-t \right) ^{\frac{1}{n}}e^{\frac{\pi i}{n}}=0\}}$$  

With Eq. (\ref{wrq_intg2}), we can express the two roots of the quadratic equation $z^2+xz-1=0$ in an integral form,

\begin{align}\label{wrq_quadratic1}
	\begin{aligned}
		z_1(x) =&\left\{1+\frac{1}{8\pi i}\int_0^1{t^{-\frac{3}{4}}}\left( 1-t \right) ^{-\frac{5}{4}}\left[ e^{\frac{\pi i}{4}}\log \left( 1-ixt^{\frac{1}{2}}\left( 1-t \right) ^{\frac{1}{2}} \right) \right.\right. \\
		& -\left. \left. e^{-\frac{\pi i}{4}}\log \left( 1+ixt^{\frac{1}{2}}\left( 1-t \right) ^{\frac{1}{2}}\right) \right] \mathrm{d}t\right\}^2 \\
	\end{aligned}
\end{align}  

\begin{align}\label{wrq_quadratic2}
	\begin{aligned}
		z_2(x)=\varepsilon z_1(\varepsilon x) =&-\left\{1+\frac{1}{8\pi i}\int_0^1{t^{-\frac{3}{4}}}\left( 1-t \right) ^{-\frac{5}{4}}\left[ e^{\frac{\pi i}{4}}\log \left( 1+ixt^{\frac{1}{2}}\left( 1-t \right) ^{\frac{1}{2}} \right) \right.\right. \\
		& -\left. \left. e^{-\frac{\pi i}{4}}\log \left( 1-ixt^{\frac{1}{2}}\left( 1-t \right) ^{\frac{1}{2}}\right) \right] \mathrm{d}t\right\}^2 \\
	\end{aligned}
\end{align}    
where $\varepsilon=\exp(\pi i)$; $x$ has a convergent domain of $\mathbb{C} ^p\backslash \left( \varSigma'' _-\cup \varSigma'' _+ \right)$, and $$\varSigma'' _-=\bigcup_{t\in \left[ 0;1 \right]}{\{ 1-ixt^{\frac{1}{2}}\left( 1-t \right) ^{\frac{1}{2}}=0 \}}$$ and $$\varSigma'' _+=\bigcup_{t\in \left[ 0;1 \right]}{\{ 1+ixt^{\frac{1}{2}}\left( 1-t \right) ^{\frac{1}{2}}=0\}}$$

It should be pointed out that only the trinomial algebraic equation $z^n+xz^{n-1}-1=0$ requires this minor modification. For other polynomial algebraic equations, the integral formula (\ref{Mikhalkin_intg1}) is correct.  

\section{Numerical examples and discussions}

Based on Eq. (\ref{wrq_intg2}), Eq. (\ref{Mikhalkin_intg1}) and Eq. (\ref{other_roots}), we can solve a large amount of polynomial algebraic equations of $n^{th}$ degree with any $n$ not only theoretically but also numerically. In this section we only give some numerical examples. Because these numerical roots are from analytical Eq. (\ref{wrq_intg2}), Eq. (\ref{Mikhalkin_intg1}) and Eq. (\ref{other_roots}), we call them semi-analytical roots. As for the pure theoretical examples, the interested readers can look up Mikhalkin (2006). 

It should be pointed out that our semi-analytical calculations here are not necessarily better than the existing numerical methods mentioned in the introduction, due to the precision and instability of numerical integration when there are singularities at the interval boundaries like ours here. We will take just one example about this, in which the {\it QuadGK} package of {\it Julia} will be used to integrate numerically $\int_0^1 t^{a-1}{\rm d}t$($=\frac{1}{a}$) where $0<a<1$. It is known that {\it QuadGK} implements an adaptive Gauss-Kronrod procedure and is accurate. When we input "${\rm quadgk}(x->x^\wedge -0.9796,0,1)$" in the REPL of Julia, we get the integral value is $49.019588870393704$ with an error of $7.297369528484689\times 10^{-7}$. However, we only get a "DomainError" when we input "${\rm quadgk}(x->x^\wedge -0.9797,0,1)$". 

In the following 9 examples will be taken, together with the corresponding results from roots formula (for quadratic equation) or eigenvalues of the corresponding companion matrix of polynomials with QR algorithm.
 
Firstly, we investigate the roots of two quadratic equations: $x^2+0.5x+1=0$ and $x^2+(10.0+2.0i)x+2=0$. Table 1 show the corresponding roots from Eq. (\ref{eq2}) and Eq. (\ref{wrq_quadratic1},\ref{wrq_quadratic2}), respectively. It can be seen that the differences between these two types of roots are small. And this shows Eq. (\ref{wrq_quadratic1}) and Eq. (\ref{wrq_quadratic2}) are correct.   

\begin{table}
	\caption{Comparison of roots for two quadratic equations}
	\begin{tabular}{ccc}
		\hline
		Equations& Roots from Eq. (\ref{eq2})& Roots from Eq. (\ref{wrq_quadratic1},\ref{wrq_quadratic2})\\
		\hline
		{\footnotesize $x^2+0.5x+1=0$} &{\footnotesize -0.250000000000000 + 0.968245836551854i} &{\footnotesize -0.249984495758362 + 0.968247805793722i} \\
		   \           &{\footnotesize -0.250000000000000 -- 0.968245836551854i} &{\footnotesize -0.249984495758362 -- 0.968247805793722i}  \\
		  {\footnotesize $x^2+(10.0+2.0i)x+2=0$} &{\footnotesize -0.195518136823225 + 0.0406949474243507i} &{\footnotesize -0.195440670967500 + 0.0407871031132125i} \\
		   \           &{\footnotesize -9.80448186317678 -- 2.04069494742435i} &{\footnotesize -9.80377786741337 -- 2.04107071213056i}  \\
		   \hline
	\end{tabular}	
\end{table}

Further, we investigate the following polynomial equations:  $z^{23}+0.5iz^{22}-1=0$, $z^{1000}+(0.5-0.37i)z^{999}-1=0$, $z^3+6z+2=0$,
$z^4-6z^2-24z+16=0$, $z^5+0.5iz^4-6iz^3-2.4z^2+z+6i=0$, $z^7+0.3z^5+0.5iz^4-6iz^3-2.4z^2+z+5=0$, and $z^{200}+0.5z^{199}+0.8z^{99}-1=0$, respectively. Fig. \ref{fig_z23}-\ref{fig_z200} show the roots of these equations. Except that the roots of $z^{23}+0.5iz^{22}-1=0$ and $z^{1000}+(0.5-0.37i)z^{999}-1=0$ are from Eq. (\ref{wrq_intg2}), roots for all other equations are from Eq. (\ref{Mikhalkin_intg1}). For comparison, we also show the roots from eigenvalues of the corresponding companion matrix. It can be seen that these two classes of roots are consistent well, except those in Fig. \ref{fig_z4} and Fig. \ref{fig_z5}. In Fig. \ref{fig_z4}, there is a root from Eq. (\ref{Mikhalkin_intg1}) is different from the eigenvalue of corresponding companion matrix. In Fig. \ref{fig_z5}, all the roots from Eq. (\ref{Mikhalkin_intg1}) are slightly different from those eigenvalues of the corresponding companion matrix, and these differences are systematic which can also be seen in other Figures except Fig. \ref{fig_z4}. We attribute these differences to the problems in numerical integration, even if the Tanh-Sinh quadrature is used (Takahasi and Mori, 1974). 

These 9 examples show that Eq. (\ref{wrq_intg2}), Eq. (\ref{Mikhalkin_intg1}) and Eq. (\ref{other_roots}) can be used directly to calculate the roots of polynomial algebraic equation of $n^{th}$ degree with any $n$, if the related numerical integration is correct. If not so, we can also use the roots from these formulae as the initial values for other numerical method mentioned in the introduction.

It can be also seen that Eq. (\ref{wrq_intg2}), Eq. (\ref{Mikhalkin_intg1}) and Eq. (\ref{other_roots}) do not converge at some points (Those in $\left( \varSigma _-\cup \varSigma _+ \right)$) on the hypercomplex plane. If the coefficients of polynomial equations consist of these points, then they cannot be solved by Eq. (\ref{wrq_intg2}), Eq. (\ref{Mikhalkin_intg1}) and Eq. (\ref{other_roots}) and another methods, such as analytical continuation, are required. As for the analytical continuation from series (\ref{Mellin1921_series1}), the interested readers are referred to Antipova and Mikhalkin (2012) and references therein.

\section{Conclusions}

A closed solution (Eq. (\ref{wrq_intg2})) for the roots of polynomial trinomial algebraic equation $z^n+xz^{n-1}-1=0$ is obtained. This solution, together with Mikhalkin's integral formula (Eq. (\ref{Mikhalkin_intg1})) and Eq. (\ref{other_roots}), not only provides a relatively simple analytical expression for the solution of a large amount of polynomial algebraic equations of $n^{th}$ degree with any $n$, but also provides a semi-analytical solution to these algebraic equations when the integral involved exists and is obtained exactly. In addition, the semi-analytical solution can provide initial values for other numerical methods.

\vspace{5em}

%\section{}
%\subsection{}

\ \ 

\end{document}